# Mass gap for the fractional $\Phi_2^4$ model

P. Duch
Institute of Mathematics
Polish Academy of Sciences
Poland
*Email:* pduch@impan.pl

M. Gubinelli
Mathematical Institute
University of Oxford
United Kingdom
*Email:* gubinelli@maths.ox.ac.uk

**Abstract**

We present a new, simple and self-contained technique for proving a mass gap for Euclidean quantum field theories in the high-temperature regime using the stochastic quantisation approach. As an application, we establish exponential decay of spatial correlations for the fractional $\Phi^4$ model on the plane in the Wick-ordering regime. The key idea is to analyse the propagation of perturbations using the linearised stochastic quantisation equation. Spatial propagation is controlled by an auxiliary random walk, while the potentially destabilising effect of the renormalised nonlinearity is controlled by combining dissipative energy estimates for the nonlinear dynamic with the Boué–Dupuis variational formula.

## Table of contents

## 1 Introduction and main results

Set $\Lambda_\varepsilon = \varepsilon \mathbb{Z}^2$ and $\Lambda_{\varepsilon,\ell} = \Lambda_\varepsilon / \ell \mathbb{Z}^2$, where the ultraviolet and volume cutoffs are

$$\varepsilon = 2^{-n}, \qquad \ell = ML, \qquad n \in \mathbb{N}_0, \quad M \in \mathbb{N}. \tag{1.1}$$

Here $L \in \mathbb{N}$ is a fixed constant chosen below. For $z \in \Lambda_\varepsilon$, we define the nearest neighbour Laplacian

$$\Delta_\varepsilon f(z) := \varepsilon^{-2} \sum_{j=1}^{2} (f(z + \varepsilon e_j) + f(z - \varepsilon e_j) - 2 f(z)).$$

For $\sigma\in(0,2]$, the operator $\mathcal{A}_\varepsilon:=(I-\Delta_\varepsilon)^{\sigma/2}$ is defined on $\Lambda_\varepsilon$ as the Fourier multiplier. Write $\nu_{\varepsilon,\ell}$ for the Gaussian measure on $\mathbb{R}^{\Lambda_{\varepsilon,\ell}}$ with covariance $\mathcal{A}_\varepsilon^{-1}$. The interacting fractional $\Phi_2^4$ measure on $\mathbb{R}^{\Lambda_{\varepsilon,\ell}}$ is

$$\mu_{\varepsilon,\ell}(\mathrm{d}\phi):=\frac{1}{\mathcal{Z}_{\varepsilon,\ell}}\exp\left(-\frac{\lambda}{4}\varepsilon^2\sum_{z\in\Lambda_{\varepsilon,\ell}}(\phi(z)^4-6\,C_{\varepsilon,\ell}\,\phi(z)^2+3\,C_{\varepsilon,\ell}^2)\right)\nu_{\varepsilon,\ell}(\mathrm{d}\phi). \tag{1.2}$$

We study the above measure in the regime in which the so-called Wick renormalisation is sufficient to control the continuum limit $\varepsilon\searrow 0$. Hence, the counterterm coincides with

$$C_{\varepsilon,\ell}:=\sum_{k\in\mathbb{Z}^2}G_\varepsilon(k\,\ell), \tag{1.3}$$

where $G_\varepsilon$ is the kernel of $\mathcal{A}_\varepsilon^{-1}$ on $\Lambda_\varepsilon$. Let $A\subset\Lambda_{\varepsilon,\ell}$. An observable $F=F_{\varepsilon,\ell}\in C^1$ is called $A$-supported if $\operatorname{supp}\nabla F(\phi)\subseteq A$, where $\nabla F$ is the gradient relative to $\langle\cdot,\cdot\rangle$. We denote by $d_\ell$ the $\ell^\infty$ distance on the torus $\mathbb{T}_\ell^2=\mathbb{R}^2/(\ell\mathbb{Z}^2)$ and write $d_\ell(A,B)$ for the distance between the $\ell$-periodic sets $A,B\subset\mathbb{T}_\ell^2$.

**Theorem 1.1.** *Let $\sigma\in(4/3,2]$. There exists $\lambda_0>0$ such that the following holds. For every $R<\infty$ there is $C_R<\infty$ such that, for every $\lambda\in(0,\lambda_0]$, every pair $(\varepsilon,\ell)$ as in (1.1), and all bounded observables $F,G$ with bounded gradients and support sets $A,B$ of diameter at most $R$,*

$$|\mathrm{Cov}_{\mu_{\varepsilon,\ell}}(F,G)|\leqslant C_R\,e^{-d_\ell(A,B)/(2L)}\sup_{\phi\in L^2}\|\nabla F(\phi)\|_2\sup_{\phi\in L^2}\|\nabla G(\phi)\|_2. \tag{1.4}$$

The assumption that $\lambda>0$ is small enough corresponds to the high-temperature regime. The case $\sigma=2$ corresponds to the standard $\Phi_2^4$ model. The exponential decay relies on our choice of the fractional operator $(1-\Delta_\varepsilon)^{\sigma/2}$. For $\sigma<2$, the alternative natural choice $1+(-\Delta_\varepsilon)^{\sigma/2}$ gives a Gaussian measure whose covariance decays only polynomially.

The proof technique does not use the scalar nature of the model and extends to the vectorial model as briefly discussed in Appendix C. By a straightforward adaptation of the technique developed in [DGR26], after a suitable extension to $\mathscr{S}'(\mathbb{R}^2)$, the family of measures $(\mu_{\varepsilon,\ell})_{\varepsilon,\ell}$ is tight in the entire subcritical regime $\sigma\in(1,2]$, for a suitable choice of counterterm that coincides with (1.3) when $\sigma\in(4/3,2]$. The bound (1.4) holds for any accumulation point for $\sigma\in(4/3,2]$. The extension for $\sigma\in(1,4/3]$, that is beyond the Wick-ordering regime, requires additional technical arguments, which we leave for future work.

The proof of the above theorem is presented in Section 6. It uses the Helffer–Sjöstrand covariance representation [HS94] combined with the analysis of the parabolic stochastic quantisation equation of the model,

$$\mathrm{d}\Phi_t=-[\mathcal{A}_\varepsilon\,\Phi_t+\lambda(\Phi_t^3-3\,C_{\varepsilon,\ell}\,\Phi_t)]\mathrm{d}t+\sqrt{2}\,\mathrm{d}W_t,\qquad \Phi_0=\phi, \tag{1.5}$$

where $W$ is standard Brownian motion in $L^2=L^2\,(\Lambda_{\varepsilon,\ell},\varepsilon^2\text{counting measure})$ extended periodically to $\Lambda_\varepsilon$. In particular,

$$\mathbb{E}[W_t(z)W_s(y)]=\varepsilon^{-2}\,1_{\{z=y\,(\mathrm{mod}\,\ell)\}}\min(t,s).$$

At fixed cutoffs (1.5) is a globally well-posed finite-dimensional SDE. Its unique invariant measure is (1.2), and the Markov semigroup $P_t$ is reversible in $L^2(\mu_{\varepsilon,\ell})$. Write $t\mapsto\Phi_t^\phi$ for the solution of (1.5) started from $\phi$ at time 0. For an arbitrary initial time $s$, we write $t\mapsto\Phi_{s,t}^\phi$ for the solution satisfying $\Phi_{s,s}^\phi=\phi$. Its derivative with respect to the initial condition is

$$U_{s,t}^\phi=D_\phi\,\Phi_{s,t}^\phi.$$

Note that $t\mapsto U_{s,t}^\phi f$ solves

$$\partial_t v_t+[\mathcal{A}_\varepsilon+3\,\lambda\,((\Phi_{s,t}^\phi)^2-C_{\varepsilon,\ell})]v_t=0,\qquad v_s=f. \tag{1.6}$$

We divide the lattice $\Lambda_{\varepsilon,\ell}$ into blocks of size $L$ and denote by $\chi_\alpha$ a bump function supported in a block centred at a lattice point $\alpha$; see Section 2.1. Our main technical result is the following contractive estimate for the linearised dynamic.

**Theorem 1.2.** *Let $\sigma \in (4/3, 2]$. There exist $\lambda_0, C > 0$, uniform in the cutoffs and block indices, such that for every $\lambda \in (0, \lambda_0]$, every pair $(\varepsilon, \ell)$ as in (1.1), every $T \geqslant 0$ and all block indices $\alpha, \beta$, we have*

$$\mathbb{E} \sup_{\psi \in \mathfrak{S}_1} \sup_{f \in \mathfrak{B}^2(\Lambda_{\varepsilon,\ell})} \|\chi_\beta U^{\psi}_{1,1+T}(\chi_\alpha f)\|^2_{L^2(\Lambda_{\varepsilon,\ell})} \leqslant C e^{-T/4 - d_\ell(\beta,\alpha)/L}, \tag{1.7}$$

*where $\mathfrak{B}^2(\Lambda_{\varepsilon,\ell})$ is the unit ball in $L^2(\Lambda_{\varepsilon,\ell})$ centred at zero and $\mathfrak{S}_1 := \{\Phi^{\phi}_1 : \phi \in L^2(\Lambda_{\varepsilon,\ell})\}$ is the random set of fields that can be reached at time one from arbitrary initial conditions.*

**Proof.** Proposition 5.2 proves the hypothesis (2.14) of Proposition 2.2, which then implies the stated bound. □

The bound (1.7) shows that a perturbation of the initial condition decays in time and has exponentially small influence on distant spatial regions uniformly in the cutoffs. This estimate extends to the infinite-volume continuum dynamic and yields exponential synchronisation of solutions driven by the same noise. Consequently, the dynamic is ergodic and admits a unique invariant measure, with every solution converging exponentially fast to the unique stationary solution. The invariant measure is invariant under all Euclidean transformations and satisfies reflection positivity. We give a self-contained proof of the results stated as theorems. The aforementioned conclusions about the limiting dynamic and its invariant measure are not proved here, but follow easily along the lines of [DHYZ25; GH21]. The proof of synchronisation by noise is outlined in Appendix D.

**Proof idea.** We first divide the lattice into blocks and dominate the propagation of the linearised equation by a conservative random walk $\mathcal{J}$ on the block lattice, whose jump rates have uniformly bounded exponential moments (see Proposition 2.1). The resulting estimate for the localised energy contains a deterministic dissipative term and a random contribution that may weaken the dissipation, measured by the local Rayleigh quotient $\hat{\mathcal{R}}_\alpha$ introduced in (2.4). Note that $\hat{\mathcal{R}}_\alpha$ involves the Hessian of the potential of the model, which is not bounded from below because of the presence of the counterterm. This is the principal reason why estimating $\hat{\mathcal{R}}_\alpha$ and the corresponding accumulated loss of dissipation is highly nontrivial. The Feynman–Kac representation reduces the problem to estimating exponential moments of the accumulated loss of dissipation along paths of the walk (see Proposition 2.2). The key tool is the Boué–Dupuis variational formula [BD98], which allows us to estimate this exponential moment by studying how the loss of dissipation changes when an arbitrary adapted drift $u$ is added to the driving Wiener process (see Proposition 5.2). More precisely, the formula (5.4) expresses the logarithm of the exponential moment as the supremum of the expected loss of dissipation for the equation with additional drift $u$, minus a cost functional equal to the squared $L^2$-norm of $u$. Thus, it suffices to bound the expected loss of dissipation uniformly over all such drifts, allowing an additional term at most quadratic in $u$, with a sufficiently small coefficient to be absorbed by the cost functional (see Lemma 5.1). To establish this estimate, we write the solution of the stochastic quantisation equation (1.5) with the additional drift $u$ as $\Phi^u = Z + \Psi^u$, separating the explicitly constructed random field $Z$ independent of the control $u$ from the more regular remainder $\Psi^u$. The stochastic bounds of Section 3 and the dissipative energy estimates of Section 4 control the loss of dissipation uniformly in the initial condition after a positive time. The coercive bound from Proposition 4.1 depends quadratically on the drift $u$, which, as explained above, is crucial for the application in the Boué–Dupuis representation.

**Related literature.** The mass gap for scalar Euclidean field theories with polynomial interactions in two dimensions was studied in the classical works [Sim73; GJS74; GRS75]. For weakly coupled $\Phi^4_3$, the existence of the infinite-volume theory and exponential decay of correlations were established in [FO76; MS76; BFS83]. The mass gap for two-dimensional Euclidean fields with exponential interactions was proved in [AH74], while more recently a proof based on stochastic quantisation was given in [GHR25]. Uniqueness of ultraregular Gibbs states for weakly coupled $P(\phi)_2$ models was proved in [AHZ89]. Logarithmic Sobolev inequalities for the $\Phi^4_2$ and $\Phi^4_3$ measures throughout the entire high-temperature phase were established in [BD24]. Building on this result, [BDW25] proved uniqueness of the invariant measure for the infinite-volume $\Phi^4_2$ stochastic quantisation dynamic in the same regime. Ergodicity of the infinite-volume $\Phi^4_3$ stochastic quantisation dynamic at high temperature, together with exponential decay of spatial correlations of its invariant measure, was established in [DHYZ25]. For the construction of the fractional $\Phi^4_3$ model and its stochastic dynamic, we refer to [Duc25; DGR26].

**Comparison with [DHYZ25].** Our approach differs substantially from that of [DHYZ25], where the linearised evolution is estimated by a time-chopping argument. This technique was originally introduced in [CLL13] to study integrability and tail estimates for rough differential equations, and was subsequently applied to the dynamic $\Phi^4_2$ model on a torus in [KT25]. The main idea is to divide the time interval into sufficiently many random subintervals such that on every subinterval the stochastic terms are small enough to obtain a pathwise contractive estimate by fixed point arguments. The number of subintervals needed to reach a deterministic time is controlled with the use of the strong Markov property and probabilistic bounds for the stopping times. Note that in this approach the stochastic objects have to be constructed with zero initial condition on every subinterval to ensure smallness. This excludes the possibility of working with stationary stochastic data. In contrast, our approach avoids completely the stopping-time argument and instead controls the possible loss of dissipation for the dynamic described by the linearised equation directly. The key ingredients are the Boué–Dupuis variational representation combined with elementary energy estimates for the nonlinear dynamic driven by the noise shifted by the Boué–Dupuis control. Moreover, an auxiliary random walk representation is used to analyse the spatial propagation of perturbations.

**Notation.** We identify functions $f\colon \Lambda_{\varepsilon,\ell} \to \mathbb{R}$ with periodic functions on $\Lambda_\varepsilon$. Kernels of operators are defined with respect to $\varepsilon^2$ times counting measure. The convolution is defined by

$$(K_\varepsilon * f)(z) := \varepsilon^2 \sum_{y \in \Lambda_\varepsilon} K_\varepsilon(z-y)\, f(y).$$

The $L^2$ inner product and $L^p$ norms are over $\Lambda_{\varepsilon,\ell}$. In particular,

$$\langle f, g\rangle := \varepsilon^2 \sum_{z \in \Lambda_{\varepsilon,\ell}} f(z)\, g(z), \qquad \|f\|_2 = \|f\|_{L^2(\Lambda_{\varepsilon,\ell})} := \langle f, f\rangle^{1/2}.$$

We stress that sums in inner products, norms and quadratic forms always run over $\Lambda_{\varepsilon,\ell}$, whereas the sum in the definition of the convolution is over $\Lambda_\varepsilon$.

**Acknowledgements.** M.G. acknowledges support from the UKRI Frontier Research / EPSRC grant "Stochastic Analysis of Quantum fields" (EP/Z534328/1) and of the Simons Collaboration on "Probabilistic QFT". This paper has been written in TEXMACS (`www.texmacs.org`).

**AI disclosure statement.** An AI model helped with the derivation of the random-walk estimates stated in Lemmas A.1 and A.3. It was also used to review the entire manuscript after its completion, correct some mathematical arguments, improve presentation as well as prepare the bibliographic entries. The authors take full responsibility for the content of the paper.

# 2 Spatial propagation

## 2.1 Block decomposition and random walk representation

We begin by studying how a perturbation of the initial data propagates under the linearised dynamic (1.6). We divide the spatial lattice into overlapping blocks of size $L$, indexed by their centres $\alpha$. For each block, we define the local $L^2$ energy $e_\alpha$ of the linearised solution inside a block $\alpha$. The kinetic operator $\mathcal{A}_\varepsilon$ transfers energy between blocks, and controlling this transfer is essential for establishing decay of correlations. The main result shows that the transfer of energy between blocks is dominated by an auxiliary random walk with exponentially decaying jump rates. This reduces the problem of studying spatial propagation to estimating the auxiliary evolution (2.7) on the lattice $\mathscr{B}_\ell$ of block centres.

Fix a nonnegative function $\chi \in C_c^\infty(\mathbb{R}^2)$, supported in $[-1,1]^2$ and strictly positive on $[-\frac{1}{2},\frac{1}{2}]^2$. The parameter $L > 16$ introduced below is fixed sufficiently big, depending only on $\sigma$ and $\chi$, such that the bounds (2.5), (2.6), (A.3) and (A.4) below hold true. Subsequently, the constant $\lambda_0 > 0$ is fixed sufficiently small. In what follows, constants $C \in (0,\infty)$ may change from line to line, and are always independent of the cutoffs, the coupling constant $\lambda \in (0,\lambda_0]$, and block indices.

For torus size $\ell = ML$, let

$$\mathscr{B}_\ell := L\,\mathbb{Z}^2 / \ell\,\mathbb{Z}^2 \subset \Lambda_{\varepsilon,\ell}$$

be the lattice of block centres. We reserve the symbols $\alpha, \beta, \bar{\alpha}, \bar{\beta} \in \mathscr{B}_\ell$ for block labels. For $z \in \Lambda_{\varepsilon,\ell}$, set

$$\tilde{\chi}_\alpha(z) := \sum_{n \in \mathbb{Z}^2} \chi\left(\frac{z - \alpha + \ell n}{L}\right), \qquad \chi_\alpha(z) := \frac{\tilde{\chi}_\alpha(z)}{(\sum_{\beta \in \mathscr{B}_\ell} \tilde{\chi}_\beta(z)^2)^{1/2}}.$$

By our assumptions on $\chi$, the denominator is strictly positive. Thus, $\chi_\alpha$ is smooth, and

$$\sum_{\alpha \in \mathscr{B}_\ell} \chi_\alpha^2 = 1, \qquad \operatorname{supp}\chi_\alpha \subset B_\alpha := \{z: d_\ell(z, \alpha) \leqslant L\}.$$

The support assumption gives

$$\operatorname{supp}\chi_\alpha \cap \operatorname{supp}\chi_\beta \neq \varnothing \quad \Longrightarrow \quad d_\ell(\alpha, \beta) \leqslant 2L.$$

Put

$$\mathcal{H}_\alpha^{\mathrm{loc}} := \{\varphi \in L^2: \operatorname{supp}\varphi \subset B_\alpha\}, \qquad w_\alpha(z) := \left(1 + \frac{d_\ell(z, \alpha)}{L}\right)^{-2}.$$

In two dimensions, we have

$$\varepsilon^2 \sum_z w_\alpha(z)^2 \leqslant C L^2, \qquad \sup_z \sum_\alpha w_\alpha(z)^2 \leqslant C, \qquad \chi_\alpha \leqslant C w_\alpha. \tag{2.1}$$

Let $q_\varepsilon(z, y)$ be the kernel introduced in Lemma A.1 such that

$$((I - \mathcal{A}_\varepsilon)v)(z) = \varepsilon^2 \sum_y q_\varepsilon(z, y)\,(v(y) - v(z)), \qquad q_\varepsilon(z, z) = 0.$$

For any lattice field $v$ define the local energy and the localisation error

$$e_\alpha(v) := \|\chi_\alpha v\|_2^2, \qquad \mathcal{I}_\alpha(v) := \varepsilon^4 \sum_{z, y} q_\varepsilon(z, y)\,(\chi_\alpha(z) - \chi_\alpha(y))^2\, v(z)\, v(y). \tag{2.2}$$

A direct computation yields

$$\langle \chi_\alpha^2 v, \mathcal{A}_\varepsilon v \rangle = \langle \chi_\alpha v, \mathcal{A}_\varepsilon(\chi_\alpha v) \rangle - \frac{1}{2}\mathcal{I}_\alpha(v). \tag{2.3}$$

For each deterministic initial condition $\phi$, define

$$\mathcal{R}_\alpha^\phi(t) := \inf_{\substack{\varphi \in \mathcal{H}_\alpha^{\mathrm{loc}} \\ \|\varphi\|_2 = 1}} \{\langle \varphi, \mathcal{A}_\varepsilon \varphi \rangle + 3\lambda\, \langle (\Phi_t^\phi)^2 - C_{\varepsilon,\ell}, \varphi^2 \rangle\}, \qquad \hat{\mathcal{R}}_\alpha^\phi(t) = (1 - \mathcal{R}_\alpha^\phi(t))_+. \tag{2.4}$$

Note that the local Rayleigh quotient $\mathcal{R}_\alpha^\phi(t)$ describes the smallest value of the quadratic form of the random Schrödinger operator

$$\mathcal{A}_\varepsilon + 3\lambda((\Phi_t^\phi)^2 - C_{\varepsilon,\ell}),$$

appearing in the linearised equation (1.6), on functions supported in a given block $\alpha$.

**Proposition 2.1.** *There exists a conservative Markov generator $\mathcal{J}$ on $\mathscr{B}_\ell$, with nonnegative rates $j_{\alpha\beta}$, $\beta \neq \alpha$, such that, for every real lattice field $v$,*

$$\mathcal{I}_\alpha(v) \leqslant (\mathcal{J}e(v))_\alpha + e_\alpha(v), \qquad (\mathcal{J}e)_\alpha = \sum_{\beta \neq \alpha} j_{\alpha\beta}\,(e_\beta - e_\alpha), \tag{2.5}$$

*and*

$$\sup_{\ell, \alpha} \sum_{\beta \neq \alpha} j_{\alpha\beta}\, e^{d_\ell(\alpha, \beta)/L} \leqslant \frac{1}{2}. \tag{2.6}$$

*Moreover, if $v$ solves the linearised equation (1.6) along $\Phi^\phi$, then*

$$\partial_t e_\alpha \leqslant (\mathcal{J}e)_\alpha - e_\alpha + 2\hat{\mathcal{R}}_\alpha^\phi e_\alpha. \tag{2.7}$$

**Proof.** Set

$$k_\alpha(z) = \varepsilon^2 \sum_y q_\varepsilon(z, y)\,(\chi_\alpha(z) - \chi_\alpha(y))^2.$$

By symmetry and Cauchy–Schwarz,

$$\mathcal{I}_\alpha(v) \leqslant \varepsilon^2 \sum_z k_\alpha(z)\, v(z)^2. \tag{2.8}$$

The bound $|\chi_\alpha(z) - \chi_\alpha(y)| \leqslant C L^{-1} d_\ell(z,y)$ and (A.2) imply $k_\alpha(z) \leqslant C L^{-2}$. If $d_\ell(z,\alpha) > L$, then $\chi_\alpha(z) = 0$ and $\chi_\alpha(y) \neq 0$ can only hold if $d_\ell(z,y) \geqslant d_\ell(z,\alpha) - L$. Thus,

$$k_\alpha(z) \leqslant \frac{C}{L^2} e^{-(d_\ell(z,\alpha)-L)_+/4} \leqslant \frac{C}{L^2} e^{-2d_\ell(z,\alpha)/L},$$

where the last inequality holds for $L \geqslant 8$, after increasing $C$.

If $\chi_\beta(z) \neq 0$, then $d_\ell(\alpha,\beta) \leqslant d_\ell(z,\alpha) + L$ by the triangle inequality. Using $\sum_\beta \chi_\beta^2 = 1$, we obtain

$$k_\alpha(z) \leqslant \sum_\beta a_{\alpha\beta}\, \chi_\beta(z)^2, \qquad a_{\alpha\beta} := \frac{C}{L^2} e^{-2d_\ell(\alpha,\beta)/L}.$$

It follows from (2.8) that

$$\mathcal{I}_\alpha(v) \leqslant \sum_\beta a_{\alpha\beta}\, e_\beta(v). \tag{2.9}$$

We have

$$\sup_\alpha \sum_\beta a_{\alpha\beta} + \sup_\alpha \sum_\beta a_{\alpha\beta}\, e^{d_\ell(\alpha,\beta)/L} \leqslant \frac{C}{L^2}. \tag{2.10}$$

Let $\mathcal{J}$ be a conservative Markov generator with rates $j_{\alpha\beta} = a_{\alpha\beta}$ for $\beta \neq \alpha$. Then

$$\mathcal{I}_\alpha(v) \leqslant \sum_\beta a_{\alpha\beta}\, e_\beta(v) = (\mathcal{J} e(v))_\alpha + \left(\sum_\beta a_{\alpha\beta}\right) e_\alpha(v).$$

Choosing $L$ sufficiently large and using (2.10) we obtain (2.5) and (2.6).

Finally, differentiating $e_\alpha(v)$ in time and using (1.6), the localisation identity (2.3), the definition of the Rayleigh quotient (2.4), and (2.5) gives (2.7). □

## 2.2 An abstract propagation criterion

The local energy inequality (2.7) contains two competing terms: a dissipative term $-e_\alpha$, which produces exponential decay, and a random potential $\hat{\mathcal{R}}^\phi_\alpha$, which may produce a local instability. The Feynman–Kac representation allows us to express the resulting energy bound as an average over paths of an auxiliary random walk $\Pi$ on spatial blocks. For $T \geqslant 0$, $k \geqslant 1$, and a piecewise-constant path $\pi\colon [1, 1+T] \to \mathscr{B}_\ell$, define

$$I_k^T := [k, k+1) \cap [1, 1+T), \qquad \mathsf{V}_k(\pi) := \{\pi(r)\colon r \in I_k^T\},$$

and

$$\mathcal{N}(\pi) := \sum_{k \geqslant 1} |\mathsf{V}_k(\pi)|, \qquad \mathcal{T}(\pi) := \sum_{k \geqslant 1} e^{-k/2} |\mathsf{V}_k(\pi)|. \tag{2.11}$$

Along each path, the random potential in (1.6) produces an exponential factor determined by the accumulated loss of dissipation. To preserve exponential decay, we formulate the criterion (2.14), which allows a small cost $\mathcal{N}(\pi)/8$ for each block visited by the walk during a unit time interval. It also permits an additional, potentially big, contribution $\mathcal{T}(\pi)$ arising from the coming-down-from-infinity estimate for the nonlinear dynamic (see the first term on the RHS of (4.4)). The distinction between these two contributions matters. Since the jump rates of the block walk are uniformly bounded, the exponential moments of $\mathcal{N}(\pi)$ grow at most exponentially in $T$. Thus, their growth must be sufficiently small to be dominated by the deterministic dissipation. In contrast, the exponentially decaying weights in $\mathcal{T}(\pi)$ imply that its exponential moments are bounded uniformly in $T$, even when multiplied by an arbitrarily large fixed constant. Thus the contribution involving $\mathcal{T}(\pi)$ changes only the multiplicative constant in the final estimate, whereas the contribution involving $\mathcal{N}(\pi)$ reduces the exponential decay rate. Note also that the naive bounds

$$\mathcal{R}^\phi_\alpha(t) \geqslant 1 - 3\lambda C_{\varepsilon,\ell}, \qquad \hat{\mathcal{R}}^\phi_\alpha(t) \leqslant 3\lambda C_{\varepsilon,\ell},$$

cannot be used directly to verify our criterion as $C_{\varepsilon,\ell}\to\infty$ as $\varepsilon\searrow 0$. To show (2.14) one has to exhibit stochastic cancellations in the local Rayleigh quotient $\mathcal{R}_\alpha^\phi(t)$. These cancellations will be analysed in subsequent sections.

We also introduce the path functional

$$\mathcal{O}^\phi(\pi) := \int_1^{1+T} \hat{\mathcal{R}}_{\pi(r)}^\phi(r)\,\mathrm{d}r, \qquad \mathcal{O}_*(\pi) = \sup_{\phi\in L^2} \mathcal{O}^\phi(\pi), \tag{2.12}$$

measuring the accumulated loss of dissipation, where the supremum in the definition of $\mathcal{O}_*(\pi)$ runs over deterministic initial conditions $\phi$. Since $\mathcal{O}^\phi(\pi)$ depends continuously on $\phi$ and the finite dimensional state space is separable, $\mathcal{O}_*(\pi)$ is measurable. The terminal time $T+1$ is determined by the domain of $\pi$. Let $\Pi$ be the continuous-time random walk on $\mathscr{B}_\ell$ generated by $\mathcal{J}$, independent of the noise, and denote by $\boldsymbol{E}_\beta^\Pi$ the expectation conditioned on $\Pi_0=\beta$. Define the uniform exponential jump cost

$$\kappa(b) := \sup_{\ell,\alpha} \sum_{\beta\neq\alpha} j_{\alpha\beta}\,(e^{d_\ell(\alpha,\beta)/L+b}-1), \tag{2.13}$$

where $j_{\alpha\beta}$ is the jump kernel of $\mathcal{J}$ from (2.5).

**Proposition 2.2.** *Suppose that, for every $T\geqslant 0$, and every deterministic piecewise-constant path $\pi\colon[1,1+T]\to\mathscr{B}_\ell$, one has*

$$\mathbb{E}\exp(2\,\mathcal{O}_*(\pi)) \leqslant \exp\big(C\,\mathcal{T}(\pi)+\tfrac{1}{8}\mathcal{N}(\pi)\big) \tag{2.14}$$

*with a deterministic $C$ independent of $T,\pi$ and the cutoffs. Then for every $T\geqslant 0$, (1.7) holds with deterministic $C$ independent of $\alpha,\beta,T,\varepsilon,\ell$.*

**Proof.** Let $v^\phi$ solve the linearised equation (1.6) along $\Phi^\phi$ on $[1,1+T]$, with $v_1^\phi=\chi_\alpha f$, and put $e_\alpha^\phi(t)=\|\chi_\alpha v_t^\phi\|_2^2$. The positive-semigroup comparison applied to (2.7) and backward Feynman–Kac formula yield

$$e_\beta^\phi(1+T) \leqslant e^{-T}\boldsymbol{E}_\beta^\Pi\bigg[\exp\bigg(2\int_0^T \hat{\mathcal{R}}_{\Pi_r}^\phi(1+T-r)\,\mathrm{d}r\bigg)\,e_{\Pi_T}^\phi(1)\bigg]. \tag{2.15}$$

The inequality holds pathwise for every initial condition. Since $v_1^\phi=\chi_\alpha f$, we have

$$e_\gamma^\phi(1)=\|\chi_\gamma\chi_\alpha f\|_2^2 \leqslant e^{2-d_\ell(\gamma,\alpha)/L}\|f\|_2^2.$$

Define the random backward path $\pi(r)=\Pi_{1+T-r}$. Taking the supremum over initial conditions in (2.15) gives

$$\begin{aligned}\sup_{\phi\in L^2} e_\beta^\phi(1+T) &\leqslant e^{-T}\boldsymbol{E}_\beta^\Pi[\exp(2\,\mathcal{O}_*(\pi))\,\|\chi_{\Pi_T}\chi_\alpha f\|_2^2]\\ &\leqslant e^{2-T}\|f\|_2^2\,\boldsymbol{E}_\beta^\Pi\,[\exp(2\,\mathcal{O}_*(\pi))\,e^{-d_\ell(\Pi_T,\alpha)/L}].\end{aligned}$$

Taking the supremum over $\|f\|_2\leqslant 1$ therefore yields

$$\sup_{\phi\in L^2}\Big\|f\mapsto\chi_\beta U_{1,1+T}^{\Phi_1^\phi}(\chi_\alpha f)\Big\|_{L^2\to L^2}^2 \leqslant e^{2-T}\boldsymbol{E}_\beta^\Pi\,[\exp(2\,\mathcal{O}_*(\pi))\,e^{-d_\ell(\Pi_T,\alpha)/L}].$$

Hence by (2.14), we obtain

$$\mathbb{E}\bigg[\sup_{\phi\in L^2}\Big\|f\mapsto\chi_\beta U_{1,1+T}^{\Phi_1^\phi}(\chi_\alpha f)\Big\|_{L^2\to L^2}^2\bigg] \leqslant e^{2-T}\boldsymbol{E}_\beta^\Pi\big[\exp\big(C\,\mathcal{T}(\pi)+\tfrac{1}{8}\mathcal{N}(\pi)\big)\,e^{-d_\ell(\Pi_T,\alpha)/L}\big].$$

Lemma A.3, applied with $a=C$ and $b=1/8$, now gives

$$\mathbb{E}\bigg[\sup_{\phi\in L^2}\Big\|f\mapsto\chi_\beta U_{1,1+T}^{\Phi_1^\phi}(\chi_\alpha f)\Big\|_{L^2\to L^2}^2\bigg] \leqslant C\,e^{-d_\ell(\beta,\alpha)/L}\,e^{-[1-1/8-\kappa(1/8)]T}.$$

Finally, by (2.6),

$$\frac{1}{8}+\kappa(1/8) \leqslant \frac{1}{8}+\frac{1}{2}e^{1/8} < \frac{3}{4}. \tag{2.16}$$

This proves (1.7). □

# 3 Stochastic enhancement and remainder dynamic

The energy estimates discussed in Section 4 will be applied to a decomposition of the stochastic quantisation solution into an explicit singular stochastic part and a more regular remainder. We first introduce the stationary Gaussian field $X$ and its Wick powers. For $\sigma > 3/2$, the Wick cube $X^{:3:}$ is sufficiently regular to be treated directly as a forcing term in the energy estimates, and one may simply use the decomposition $\Phi = X + \bar{\Psi}$. To cover the entire range $\sigma > 4/3$, we introduce the second order correction $Y$ to remove $X^{:3:}$ from the remainder equation.

Let $W$ be the two sided Wiener process. Define stationary stochastic convolutions

$$(\partial_t + \mathcal{A}_\varepsilon)X = \sqrt{2}\ \dot{W}, \qquad X^{:2:} = X^2 - C_{\varepsilon,\ell}, \qquad X^{:3:} = X^3 - 3\,C_{\varepsilon,\ell}\,X, \qquad (\partial_t + \mathcal{A}_\varepsilon)Y = X^{:3:},$$

where $C_{\varepsilon,\ell}$ is the counterterm introduced in (1.3). Note that $C_{\varepsilon,\ell} = \mathbb{E}[X_t(z)^2]$. Let $p_t^\varepsilon$ be the infinite-lattice convolution kernel of $e^{-t\mathcal{A}_\varepsilon}$. We have

$$X(t,z) = \sqrt{2}\int_{-\infty}^t (p_{t-r}^\varepsilon * \mathrm{d}W_r)(z), \qquad Y(t,z) = \int_{-\infty}^t (p_{t-r}^\varepsilon * X^{:3:}(r))(z)\,\mathrm{d}r.$$

Put $Z := X - \lambda Y$ and

$$V := X^{:2:} - 2\,\lambda\,X\,Y + \lambda^2\,Y^2 = Z^2 - C_{\varepsilon,\ell}, \qquad H := -3\,\lambda\,X^{:2:}\,Y + 3\,\lambda^2\,X\,Y^2 - \lambda^3\,Y^3.$$

The fields $V$ and $H$ play the role of the random potential and the forcing term in the equation for the remainder $\Psi := \Phi - Z$.

For a deterministic initial condition $\phi$ and a control $u$, let $\Phi^{\phi,u}$ denote the solution of the stochastic quantisation equation with the additional drift $\sqrt{2}u$ and the initial condition $\Phi_0^{\phi,u} = \phi$. Define $\Psi^{\phi,u} := \Phi^{\phi,u} - Z$. We suppress the dependence on the initial condition in the notation when $\phi$ is fixed, and write $\Phi^u = \Phi^{\phi,u}$ and $\Psi^u = \Psi^{\phi,u}$.

**Lemma 3.1.** *The remainder $\Psi^u$ satisfies*

$$(\partial_t + \mathcal{A}_\varepsilon)\Psi^u + \lambda(\Psi^u)^3 = -3\,\lambda\,Z(\Psi^u)^2 - 3\,\lambda\,V\,\Psi^u - \lambda\,H + \sqrt{2}\,u. \tag{3.1}$$

**Proof.** Expand $(Z + \Psi^u)^3 - 3\,C_{\varepsilon,\ell}\,(Z + \Psi^u)$ using $Z = X - \lambda\,Y$, and cancel the cubic Wick term using the equation for $Y$. □

Fix $0 < \delta < (3\sigma - 4)/12$ and set

$$\eta = 1 - \frac{\sigma}{2} + \delta. \tag{3.2}$$

Then $\sigma - 4\eta > 0$ and $2\eta < \sigma/2$. Choose

$$p \in \mathbb{N}, \qquad p \geqslant \max\left\{8, \frac{2\sigma}{3\sigma - 4 - 6\delta}\right\}, \qquad \theta = \frac{1}{10p}$$

and define

$$\omega_\alpha(z) = \left(1 + \frac{d_\ell(z,\alpha)}{L}\right)^{-\theta}.$$

The corresponding weighted Besov norms are defined in Appendix B. The constant $p$ plays the role of the exponent of the size of the enhancement $\mathfrak{A}_\alpha(t)$ introduced below, and is chosen sufficiently large to absorb the powers of $\mathfrak{A}_\alpha(t)$ arising later in the energy and Rayleigh-quotient estimates in Section 4. The choice of $\theta$ ensures $8\,\theta < 2$, which guarantees the spatial integrability of $w_\alpha^2\,\omega_\alpha^{-8}$ used in Lemma B.5.

**Definition 3.2.** *For $\alpha \in \mathscr{B}_\ell$ and $t \in \mathbb{R}$, define*

$$\begin{aligned}\mathfrak{A}_\alpha(t) :=\ & 1 + \|X(t)\|_{\mathcal{C}^{-\eta}(\omega_\alpha)} + \|X^{:2:}(t)\|_{\mathcal{C}^{-2\eta}(\omega_\alpha)}\\ & + \|Y(t)\|_{\mathcal{C}^{\sigma-3\eta}(\omega_\alpha)} + \|X^{:2:}\,Y(t)\|_{\mathcal{C}^{-2\eta}(\omega_\alpha)}.\end{aligned}$$

Note that $\mathfrak{A}_\alpha(t)$ does not depend on $\lambda$. By Lemma B.3, there is a deterministic constant $C$ such that for all $\alpha \in \mathscr{B}_\ell$, $t \in \mathbb{R}$, and $0 < \lambda \leqslant 1$,

$$\|Z(t)\|_{\mathcal{C}^{-\eta}(\omega_\alpha)} \vee \|V(t)\|_{\mathcal{C}^{-2\eta}(\omega_\alpha^2)}^{1/2} \vee \|H(t)\|_{\mathcal{C}^{-2\eta}(\omega_\alpha^3)}^{1/3} \leqslant C\,\mathfrak{A}_\alpha(t). \tag{3.3}$$

**Proposition 3.3.** *For every $p \geqslant 1$, there is $C_p < \infty$ such that*

$$\sup_{\varepsilon, \ell} \sup_{\alpha \in \mathscr{B}_\ell} \sup_{t \in \mathbb{R}} \mathbb{E}\, \mathfrak{A}_\alpha(t)^p \leqslant C_p. \tag{3.4}$$

**Proof.** The result can be proved by adapting the proof technique from Section 6 of [ZZ18], where analogous stochastic estimates are proved in the more complicated setting of the $\Phi_3^4$ model. □

# 4 Pathwise energy and Rayleigh estimates

## 4.1 Weighted energy inequality

To establish the exponential-moment bound (2.14) required by Proposition 2.2, we shall apply in Section 5 the Boué–Dupuis variational formula. In this section, we derive suitable pathwise energy estimates for solutions of the stochastic quantisation equation (1.5) with an additional drift $u$, which will play the role of the control in the Boué–Dupuis representation.

For $q \geqslant 1$ define

$$\|f\|_q^q = \varepsilon^2 \sum_z |f(z)|^q, \qquad \|f\|_{H^a} = \|(I - \Delta_\varepsilon)^{a/2} f\|_2, \qquad a \geqslant 0.$$

Let $u \in L^2([0, t]; L^2)$ be a control. Recall that the controlled dynamic satisfies (1.5) with an additional drift $\sqrt{2}\, u$ and the initial condition $\Phi_0^u = \phi$. Define

$$\mathcal{E}_\alpha^u(t) := \|w_\alpha \Psi^u(t)\|_2^2, \qquad \mathcal{D}_\alpha^u(t) := \|w_\alpha \Psi^u(t)\|_{H^{\sigma/2}}^2 + \lambda \|w_\alpha^{1/2} \Psi^u(t)\|_4^4, \tag{4.1}$$

where $\Psi^u = \Phi^u - Z$.

**Proposition 4.1.** *For $\lambda \in (0, \lambda_0]$, there is $C > 0$ uniform in $\varepsilon, \ell, \alpha$, the deterministic initial condition $\phi$, and the control $u$, such that*

$$\frac{\mathrm{d}}{\mathrm{d}t} \mathcal{E}_\alpha^u + \frac{1}{2} \mathcal{E}_\alpha^u + \frac{\lambda}{2 L^2} (\mathcal{E}_\alpha^u)^2 + \frac{1}{2} \mathcal{D}_\alpha^u \leqslant C \|w_\alpha u\|_2^2 + C \lambda\, \mathfrak{A}_\alpha(t)^p. \tag{4.2}$$

**Proof.** Write

$$E = \mathcal{E}_\alpha^u, \qquad K = \|w_\alpha \Psi^u\|_{H^{\sigma/2}}^2, \qquad P = \|w_\alpha^{1/2} \Psi^u\|_4^4.$$

Testing (3.1) against $2 w_\alpha^2 \Psi^u$ and applying Lemmas A.2 and B.5 gives

$$\frac{\mathrm{d}}{\mathrm{d}t} E + \frac{7}{4} K + \frac{7}{4} \lambda P \leqslant \frac{1}{8} E + 2\sqrt{2} |\langle w_\alpha \Psi^u, w_\alpha u \rangle| + C \lambda\, \mathfrak{A}_\alpha^p.$$

Since $K \geqslant E$ and $2\sqrt{2} |\langle w_\alpha \Psi^u, w_\alpha u \rangle| \leqslant \frac{1}{4} E + C \|w_\alpha u\|_2^2$, we obtain

$$\frac{\mathrm{d}}{\mathrm{d}t} E + \frac{11}{8} K + \frac{7}{4} \lambda P \leqslant C \|w_\alpha u\|_2^2 + C \lambda\, \mathfrak{A}_\alpha^p.$$

Finally, (A.4) gives $L^{-2} E^2 \leqslant 2 P$. Hence

$$\frac{1}{2} E + \frac{1}{2} K \leqslant \frac{11}{8} K, \qquad \frac{\lambda}{2 L^2} E^2 + \frac{\lambda}{2} P \leqslant \frac{7}{4} \lambda P.$$

This proves (4.2). □

## 4.2 Coming down from infinity and dissipation

The energy inequality obtained above depends on the initial condition, whereas the estimate stated in Theorem 1.2 is uniform in the initial data. To obtain the desired bound, we use the coming-down-from-infinity estimate (4.7) below to bound the initial energy in terms of the forcing at earlier times, with a weight decaying exponentially in time.

**Lemma 4.2.** *Let $T\geqslant 0$, $k\geqslant 1$, and recall the intervals $I_k^T$ introduced in Section 2. Define*

$$K_k^T(r) := e^{-(k-r)/2}\,1_{\{0\leqslant r<k\}} + 1_{\{r\in I_k^T\}}. \tag{4.3}$$

*Then*

$$\lambda\int_{I_k^T}\mathcal{D}_\alpha^u(r)\,\mathrm{d}r \leqslant C\,L^2\,e^{-k/2} + C\,\lambda\int_0^{1+T} K_k^T(r)\,(\|w_\alpha\,u(r)\|_2^2 + \lambda\,\mathfrak{A}_\alpha(r)^p)\,\mathrm{d}r. \tag{4.4}$$

*Moreover,*

$$\sup_{r\in[0,1+T]}\sum_{k\geqslant 1} K_k^T(r) \leqslant C, \qquad \sup_{k\geqslant 1}\int_0^{1+T} K_k^T(r)\,\mathrm{d}r \leqslant C. \tag{4.5}$$

*The constants $C$ above are deterministic and independent of $T, k, \varepsilon, \ell, \alpha$, the initial condition $\phi$, the control $u$, and $\lambda\in(0,\lambda_0]$.*

**Proof.** Let $C$ be as in (4.2) and put

$$h(r) = C\int_0^r e^{-(r-s)/2}\,(\|w_\alpha\,u(s)\|_2^2 + \lambda\,\mathfrak{A}_\alpha(s)^p)\,\mathrm{d}s.$$

Since $\partial_r h + h/2$ coincides with the forcing in (4.2) and $h\geqslant 0$, setting $g = (\mathcal{E}_\alpha^u - h)_+$, we arrive at

$$\partial_r g + \frac{1}{2}\,g + \frac{\lambda}{2L^2}\,g^2 \leqslant 0, \tag{4.6}$$

where we used $\mathcal{E}_\alpha^u \geqslant g$. If $g>0$, the inverse $1/g$ satisfies

$$\partial_r\frac{1}{g} \geqslant \frac{1}{2g} + \frac{\lambda}{2L^2}.$$

Integrating this inequality gives

$$\frac{1}{g(r)} \geqslant \frac{e^{r/2}}{g(0)} + \frac{\lambda}{L^2}\,(e^{r/2}-1) \geqslant \frac{\lambda}{L^2}\,(e^{r/2}-1).$$

If $g$ vanishes at some point in time, it must be identically zero by (4.6). Consequently,

$$(\mathcal{E}_\alpha^u(r) - h(r))_+ \leqslant \frac{L^2}{\lambda\,(e^{r/2}-1)}.$$

Thus, independently of the initial energy,

$$\mathcal{E}_\alpha^u(r) \leqslant \frac{C\,L^2}{\lambda\,(e^{r/2}-1)} + C\int_0^r e^{-(r-s)/2}\,(\|w_\alpha\,u(s)\|_2^2 + \lambda\,\mathfrak{A}_\alpha(s)^p)\,\mathrm{d}s, \qquad r>0. \tag{4.7}$$

The bound (4.4) follows by integrating (4.2) over $I_k^T$, discarding the nonnegative terms, and bounding the energy at $k$ by (4.7), using $(e^{k/2}-1)^{-1}\leqslant C\,e^{-k/2}$ for $k\geqslant 1$. The kernel bounds follow by elementary estimates. □

## 4.3 Local Rayleigh-quotient bound

For a reference solution $\Phi^u = Z + \Psi^u$ starting from the initial condition $\phi$ and driven by the noise shifted by $u$, let $\mathcal{R}_\alpha^{\phi,u}$ be defined by (2.4) with $\Phi^\phi$ replaced by $\Phi^{\phi,u}$, and set $\hat{\mathcal{R}}_\alpha^{\phi,u} = (1-\mathcal{R}_\alpha^{\phi,u})_+$. The main difficulty in estimating the local Rayleigh quotient is the absence of a uniform lower bound on the random potential due to the presence of the counterterm. The decomposition $\Phi^u = Z + \Psi^u$ of the nonlinear solution separates this potential into a non-negative term quadratic in the remainder and two potentially negative terms involving the stochastic coefficients,

$$(\Phi^u)^2 - C_{\varepsilon,\ell} = (\Psi^u)^2 + V + 2\,Z\,\Psi^u. \tag{4.8}$$

We estimate the last two terms using the weighted Besov estimates, which allow us to bound their contributions by the kinetic term and $\mathcal{D}_\alpha^u$. As a result, we can bound $\hat{\mathcal{R}}_\alpha^{\phi,u}(t)$ entirely in terms of the dissipation $\mathcal{D}_\alpha^u$ and the stochastic enhancement.

**Proposition 4.3.** *For* $\lambda \in (0, \lambda_0]$ *and any* $\varphi \in \mathcal{H}_\alpha^{\rm loc}$ *such that* $\|\varphi\|_2 = 1$*, one has*

$$3\,\lambda\,\langle V(t) + 2\,Z(t)\,\Psi^u(t), \varphi^2\rangle \geqslant -\tfrac{1}{2}\langle \varphi, (\mathcal{A}_\varepsilon - I)\,\varphi\rangle - C\,(\lambda\,\mathcal{D}_\alpha^u(t) + \lambda\,\mathfrak{A}_\alpha(t)^p). \tag{4.9}$$

*In particular,*

$$\hat{\mathcal{R}}_\alpha^{\phi,u}(t) \leqslant C\,\lambda\,(\mathcal{D}_\alpha^u(t) + \mathfrak{A}_\alpha(t)^p). \tag{4.10}$$

*The constants are uniform in* $\varepsilon, \ell, \alpha, \phi, u$.

**Proof.** By (3.3) and Lemma B.4, we have

$$|\langle V, \varphi^2\rangle| \leqslant C\,\mathfrak{A}_\alpha^2\,\|\varphi\|_2\|\varphi\|_{H^{2\eta+2\delta}},$$
$$|\langle Z\,\Psi^u, \varphi^2\rangle| \leqslant C\,\mathfrak{A}_\alpha\,\|w_\alpha \Psi^u\|_{H^{\sigma/2}}\|\varphi\|_{H^{\eta+\delta/2}}^2.$$

Since $\|\varphi\|_2 = 1$, interpolation between Sobolev spaces gives

$$\|\varphi\|_{H^{2\eta+2\delta}} \leqslant \|\varphi\|_{H^{\sigma/2}}^{\nu}, \qquad \|\varphi\|_{H^{\eta+\delta/2}}^2 \leqslant \|\varphi\|_{H^{\sigma/2}}^{\gamma},$$

with

$$\nu = \frac{4\,\eta + 4\,\delta}{\sigma} < 1, \qquad \gamma = \frac{4\,\eta + 2\,\delta}{\sigma} < 1.$$

Hence, by Young's inequality

$$3\,\lambda|\langle V, \varphi^2\rangle| \leqslant C\,\lambda\,\mathfrak{A}_\alpha^2\|\varphi\|_{H^{\sigma/2}}^{\nu} \leqslant \lambda\|\varphi\|_{H^{\sigma/2}}^2 + C\,\lambda\,\mathfrak{A}_\alpha^{4/(2-\nu)},$$
$$6\,\lambda|\langle Z\,\Psi^u, \varphi^2\rangle| \leqslant C\,\lambda\,\mathfrak{A}_\alpha\,\|w_\alpha\,\Psi^u\|_{H^{\sigma/2}}\|\varphi\|_{H^{\sigma/2}}^{\gamma} \leqslant \lambda\|\varphi\,\|_{H^{\sigma/2}}^2 + \lambda\|\,w_\alpha\,\Psi^u\|_{H^{\sigma/2}}^2 + C\,\lambda\,\mathfrak{A}_\alpha^{2/(1-\gamma)}.$$

By the assumptions on $\delta$ and $p$ from Section 3,

$$\frac{4}{2-\nu} < 4 \leqslant p, \qquad \frac{2}{1-\gamma} = \frac{2\,\sigma}{3\,\sigma - 4 - 6\,\delta} \leqslant p.$$

Since $\mathfrak{A}_\alpha \geqslant 1$, $\|\varphi\|_{H^{\sigma/2}}^2 \lesssim 1 + \langle\varphi, (\mathcal{A}_\varepsilon - I)\,\varphi\rangle$ and $\|w_\alpha\,\Psi^u\|_{H^{\sigma/2}}^2 \leqslant \mathcal{D}_\alpha^u$, the above estimates imply (4.9). Finally, using (2.4) and discarding the nonnegative first term in (4.8), we obtain (4.10). □

# 5 Path estimates and Boué–Dupuis

## 5.1 Bounds along block paths

For a path $\pi: [1, 1+T] \to \mathscr{B}_\ell$ and a control $u$, define $\mathcal{O}^{\phi,u}(\pi)$ and $\mathcal{O}_*^u(\pi)$ by replacing $\hat{\mathcal{R}}_\alpha^\phi$ with $\hat{\mathcal{R}}_\alpha^{\phi,u}$ in (2.12). Explicitly, $\mathcal{O}_*^u(\pi) := \sup_{\phi\in L^2}\mathcal{O}^{\phi,u}(\pi)$. In this section, we estimate the accumulated loss of dissipation $\mathcal{O}_*^u(\pi)$ along a fixed path $\pi$ by combining the local Rayleigh-quotient bound with the energy estimates for the nonlinear dynamic. The estimate is uniform in the initial condition and the cutoffs and provides the main input for the application of the Boué–Dupuis formula in the next subsection.

**Lemma 5.1.** *For every* $\lambda \in (0, \lambda_0]$*, every* $T \geqslant 0$*, every piecewise-constant path* $\pi: [1, 1+T] \to \mathscr{B}_\ell$*, and every control* $u$*,*

$$\mathcal{O}_*^u(\pi) \leqslant C\,\mathcal{T}(\pi) + C\,\lambda\int_0^{1+T}\|u(r)\|_2^2\,\mathrm{d}r + C\,\lambda\,\mathfrak{M}(\pi), \qquad \mathbb{E}\,\mathfrak{M}(\pi) \leqslant C\,\mathcal{N}(\pi). \tag{5.1}$$

*The constant* $C$ *is deterministic and uniform in* $\varepsilon, \ell, \lambda, T, \pi, u$*. The nonnegative random variable* $\mathfrak{M}(\pi)$ *depends only on the enhancement and the path, not on* $u$ *or* $\phi$*.*

**Proof.** For $r \in I_k^T$, $\pi(r) \in \mathsf{V}_k(\pi)$. Hence (4.10) yields

$$\mathcal{O}^{\phi,u}(\pi) \leqslant C\,\lambda\sum_{k\geqslant 1}\sum_{\alpha\in\mathsf{V}_k(\pi)}\int_{I_k^T}(\mathcal{D}_\alpha^u(r) + \mathfrak{A}_\alpha(r)^p)\,\mathrm{d}r.$$

The contribution from the first term on the RHS of (4.4) sums to $C\,L^2\,\mathcal{T}(\pi)$. As $K_k^T \geqslant 1_{I_k^T}$ and $\lambda \leqslant 1$, the enhancement terms are bounded by $C\,\lambda\,\mathfrak{M}(\pi)$, where

$$\mathfrak{M}(\pi) = \sum_{k \geqslant 1} \sum_{\alpha \in \mathsf{V}_k(\pi)} \int_0^{1+T} K_k^T(r)\,\mathfrak{A}_\alpha(r)^p\,\mathrm{d}r. \tag{5.2}$$

By Proposition 3.3 and (4.5), $\mathbb{E}\,\mathfrak{M}(\pi) \leqslant C\,\mathcal{N}(\pi)$. The contribution from the term on the RHS of (4.4) involving the control $u$ is bounded independently of the number of blocks visited by the path. Indeed, by (2.1), we have

$$\sum_{\alpha \in \mathsf{V}_k(\pi)} \|w_\alpha\,u(r)\|_2^2 \leqslant C\,\|u(r)\|_2^2.$$

Summing over $k$ and using (4.5), we bound the total control contribution by $C\,\lambda\int_0^{1+T}\|u(r)\|_2^2\,\mathrm{d}r$. Since the estimates above hold pathwise with a right-hand side independent of the initial condition, taking the supremum over all initial conditions gives (5.1). □

## 5.2 Exponential path bound

We finally establish the exponential-moment estimate in the propagation criterion formulated in Proposition 2.2. This finishes the proof of the decay estimate for the linearised evolution stated in Theorem 1.2.

**Proposition 5.2.** *There exists a sufficiently small $\lambda_0 > 0$ such that, for every $\lambda \in (0, \lambda_0]$, every $T \geqslant 0$, and every piecewise-constant path $\pi\colon [1, 1+T] \to \mathscr{B}_\ell$,*

$$\mathbb{E}\exp\left(2\,\mathcal{O}_*(\pi)\right) \leqslant \exp\!\big(C\,\mathcal{T}(\pi) + \tfrac{1}{8}\mathcal{N}(\pi)\big). \tag{5.3}$$

*The expectation is over the noise, and the constant $C$ is independent of $\lambda, T, \pi$ and the cutoffs.*

**Proof.** Since $\mathcal{A}_\varepsilon \geqslant I$ and $3\lambda(\Phi^2 - C_{\varepsilon,\ell}) \geqslant -3\,\lambda\,C_{\varepsilon,\ell}$, we have $0 \leqslant \mathcal{O}^{\phi,u}(\pi) \leqslant 3\,\lambda\,C_{\varepsilon,\ell}\,T$. For deterministic initial data the functional $\mathcal{O}^\phi(\pi)$ depends only on the Wiener process after time zero. The same properties are true for $\mathcal{O}_*(\pi)$. Thus the Boué–Dupuis formula [BD98] applies to this bounded functional and yields

$$\log\mathbb{E}\,e^{2\mathcal{O}_*(\pi)} = \sup_u \mathbb{E}\left[2\mathcal{O}_*^u(\pi) - \frac{1}{2}\int_0^{1+T}\|u(r)\|_2^2\,\mathrm{d}r\right]. \tag{5.4}$$

The supremum is over controls $u$ with finite expected $L^2$ energy progressively measurable with respect to the filtration generated by the Wiener process after time zero. To bound the shifted functional, we use the decomposition $\Phi^{\phi,u} = Z + \Psi^{\phi,u}$ and use the enhancement from Section 3, which depends on the original, unshifted noise also before time zero. Note that the use of the noise before time zero plays an auxiliary role and does not affect the variational formulation. To finish the proof insert (5.1) and fix $\lambda_0$ sufficiently small that $2\,C\,\lambda_0 \leqslant 1/8$. □

# 6 Proof of the mass gap

To prove Theorem 1.1 we combine the estimate for the linearised dynamic from Theorem 1.2 with the Helffer–Sjöstrand [HS94] representation of the covariance. The estimate for the linearised evolution uses the coming-down-from-infinity property and becomes uniform in the initial condition only after a positive time. To exploit this estimate in the covariance representation below, we use stationarity to shift the tangent evolution from $[0, T]$ to $[1, 1+T]$.

**Lemma 6.1.** *Let $\Phi$ be a stationary solution of (1.5) with initial law $\mu_{\varepsilon,\ell}$. For bounded $C^1$ observables $F, G$ with bounded gradients, one has*

$$\mathrm{Cov}_{\mu_{\varepsilon,\ell}}(F, G) = \int_0^\infty \mathbb{E}_{\mu_{\varepsilon,\ell}}[\langle U_{1,1+T}^{\Phi_1}\,\nabla F(\Phi_1), \nabla G(\Phi_{1+T})\rangle]\,\mathrm{d}T. \tag{6.1}$$

**Proof.** Let $P_T$ be the Markov semigroup associated to (1.5) and write $\mu=\mu_{\varepsilon,\ell}$ for the invariant measure. Reversibility and integration by parts give

$$\frac{\mathrm{d}}{\mathrm{d}T}\int F P_T G\,\mathrm{d}\mu = -\int \langle \nabla F, \nabla P_T G\rangle\,\mathrm{d}\mu.$$

Since the finite-dimensional dynamic (1.5) is ergodic, $P_T G\to\int G\,\mathrm{d}\mu$ in $L^2(\mu)$ as $T\to\infty$. Integrating the preceding identity yields

$$\mathrm{Cov}_\mu(F,G)=\int_0^\infty\int\langle\nabla F,\nabla P_T G\rangle\,\mathrm{d}\mu\,\mathrm{d}T.$$

By the chain rule,

$$\langle\nabla P_T G(\phi), f\rangle=\mathbb{E}[\langle\nabla G(\Phi_T^\phi), U_{0,T}^\phi f\rangle].$$

Combining the last two identities and using stationarity to shift the time proves (6.1). □

**Proof of Theorem 1.1.** Let $\Phi$ be stationary with fixed-time law $\mu_{\varepsilon,\ell}$, and set

$$f=\nabla F(\Phi_1),\qquad g=\nabla G(\Phi_{1+T}).$$

By the partition of unity,

$$\langle U_{1,1+T}^{\Phi_1} f, g\rangle=\sum_{\alpha,\beta}\langle\chi_\beta U_{1,1+T}^{\Phi_1}(\chi_\alpha^2 f),\chi_\beta g\rangle.$$

Set

$$F_\alpha=\sup_{\phi\in L^2}\|\chi_\alpha\nabla F(\phi)\|_2,\qquad G_\beta=\sup_{\phi\in L^2}\|\chi_\beta\nabla G(\phi)\|_2.$$

Cauchy–Schwarz and Theorem 1.2, applied conditionally on $\Phi_0$ with initial data $\chi_\alpha f$ for the linearised dynamic, give

$$\begin{aligned}|\mathbb{E}_\mu\langle U_{1,1+T}^{\Phi_1} f, g\rangle| &\leqslant \sum_{\alpha,\beta} G_\beta\,(\mathbb{E}_\mu\|\chi_\beta U_{1,1+T}^{\Phi_1}(\chi_\alpha^2 f)\|_2^2)^{1/2}\\ &\leqslant C e^{-T/8}\sum_{\alpha,\beta}F_\alpha G_\beta e^{-d_\ell(\alpha,\beta)/(2L)}.\end{aligned}$$

As $A$ and $B$ have diameter at most $R$, there are only finitely many $\alpha$ and $\beta$ such that $\mathrm{supp}\chi_\alpha\cap A\neq\emptyset$ and $\mathrm{supp}\chi_\beta\cap B\neq\emptyset$. Moreover, for such indices,

$$d_\ell(\alpha,\beta)\geqslant d_\ell(A,B)-2L.$$

Consequently,

$$|\mathbb{E}_\mu\langle U_{1,1+T}^{\Phi_1} f, g\rangle|\leqslant C_R\,e^{-T/8-d_\ell(A,B)/(2L)}\sup_{\phi\in L^2}\|\nabla F(\phi)\|_2\sup_{\phi\in L^2}\|\nabla G(\phi)\|_2.$$

Integrating over $T\geqslant 0$ and using (6.1) proves (1.4). □

# Appendix A Random-walk estimates

In this appendix we collect estimates for the discrete random walk. We first represent the nonlocal operator $I-\mathcal{A}_\varepsilon$ as the generator of a symmetric random walk and establish exponential moment bounds for its jump kernel. These bounds are then used to control errors in the spatial localisation. Finally, we estimate exponential moments of the number of blocks visited by the random walk. The latter estimate is the key ingredient in the proof of the abstract propagation criterion (2.14).

**Lemma A.1.** *The operator $I-\mathcal{A}_\varepsilon$ is a symmetric conservative Markov generator. There is a nonnegative symmetric periodic kernel $q_\varepsilon$, with $q_\varepsilon(z,z)=0$, such that*

$$((I-\mathcal{A}_\varepsilon)f)(z)=\varepsilon^2\sum_y q_\varepsilon(z,y)\,(f(y)-f(z)).\qquad\text{(A.1)}$$

*Moreover,*

$$\sup_{\varepsilon, \ell, z} \varepsilon^2 \sum_y q_\varepsilon(z, y)\, d_\ell(z, y)^2\, e^{d_\ell(z, y)/4} \leqslant C_\sigma. \tag{A.2}$$

**Proof.** For $\sigma = 2$, the kernel is given, for $y \neq z$, by

$$q_\varepsilon(z, y) = \varepsilon^{-4} \sum_{j=1}^{2} (1_{\{y = z + \varepsilon \mathrm{e}_j \,(\mathrm{mod}\, \ell)\}} + 1_{\{y = z - \varepsilon \mathrm{e}_j \,(\mathrm{mod}\, \ell)\}}),$$

where $\mathrm{e}_1 = (1, 0)$ and $\mathrm{e}_2 = (0, 1)$. The claimed properties follow directly.

Suppose $\sigma < 2$. Let $S_s$ be the continuous-time nearest-neighbour walk on $\Lambda_\varepsilon$, started at 0, with generator $\Delta_\varepsilon$, and let $\bar{S}_s$ denote its projection onto $\Lambda_{\varepsilon, \ell}$. The identity

$$I - \mathcal{A}_\varepsilon = \frac{\sigma/2}{\Gamma(1 - \sigma/2)} \int_0^\infty e^{-s}\, s^{-1-\sigma/2}\, (e^{s\Delta_\varepsilon} - I)\, \mathrm{d}s$$

gives the desired generator representation with $q_\varepsilon(z, z) = 0$ and, for $y \neq z$,

$$q_\varepsilon(z, y) = \frac{\sigma/2}{\varepsilon^2\, \Gamma(1 - \sigma/2)} \int_0^\infty e^{-s}\, s^{-1-\sigma/2}\, \mathbb{P}(\bar{S}_s = y - z)\, \mathrm{d}s.$$

In particular, $q_\varepsilon$ is nonnegative and symmetric. Let $\varsigma \in \{-1, 1\}^2$. Since $z \mapsto e^{\vartheta \varsigma \cdot z}$ is an eigenfunction of $\Delta_\varepsilon$ with eigenvalue $4\,\varepsilon^{-2}\,(\cosh(\varepsilon\,\vartheta) - 1)$, we obtain

$$\mathbb{E}\, e^{\vartheta \varsigma \cdot S_s} = \exp\!\left(\frac{4\,s}{\varepsilon^2}\,(\cosh(\varepsilon\,\vartheta) - 1)\right).$$

Differentiating twice in $\vartheta$, summing over $\varsigma$, and using $0 < \varepsilon \leqslant 1$, we obtain

$$\mathbb{E}[|S_s|^2\, e^{|S_s|/4}] \leqslant C\,(s + s^2)\, e^{cs}, \qquad c = 4\,(\cosh(1/4) - 1) < 1.$$

Since $d_\ell(0, \bar{S}_s) \leqslant |S_s|$, it follows that

$$\begin{aligned}
&\varepsilon^2 \sum_y q_\varepsilon(z, y)\, d_\ell(z, y)^2\, e^{d_\ell(z, y)/4} \\
&\quad \leqslant \frac{\sigma/2}{\Gamma(1 - \sigma/2)} \int_0^\infty e^{-s}\, s^{-1-\sigma/2}\, \mathbb{E}\,[|S_s|^2\, e^{|S_s|/4}]\, \mathrm{d}s \\
&\quad \leqslant C_\sigma \int_0^\infty (s + s^2)\, s^{-1-\sigma/2}\, e^{-(1-c)s}\, \mathrm{d}s \leqslant C_\sigma.
\end{aligned}$$

The integral is finite because $\sigma < 2$ and $c < 1$. This proves (A.2). □

**Lemma A.2.** *For the fixed sufficiently large $L$, all $\alpha \in \mathscr{B}_\ell$, and all real lattice fields $\psi$,*

$$2\,\langle w_\alpha^2\, \psi, \mathcal{A}_\varepsilon\, \psi\rangle \geqslant 2\|w_\alpha\, \psi\|_{H^{\sigma/2}}^2 - \frac{1}{8}\,\|w_\alpha\, \psi\|_2^2 \tag{A.3}$$

*and*

$$L^{-2}\,\|w_\alpha\, \psi\|_2^4 \leqslant 2\,\|w_\alpha^{1/2}\, \psi\|_4^4. \tag{A.4}$$

**Proof.** Using the fact that $\rho \mapsto \log(1 + \rho/L)$ is Lipschitz on $[0, \infty)$ with Lipschitz constant $L^{-1}$, as well as the triangle inequality for $d_\ell$, we obtain $|\log w_\alpha(z) - \log w_\alpha(y)| \leqslant 2\, d_\ell(z, y)/L$. Consequently, for $L \geqslant 8$, we have

$$\frac{(w_\alpha(z) - w_\alpha(y))^2}{w_\alpha(z)\, w_\alpha(y)} = 2\,(\cosh(\log w_\alpha(z) - \log w_\alpha(y)) - 1) \leqslant \frac{C}{L^2}\, d_\ell(z, y)^2\, e^{d_\ell(z, y)/4}.$$

By (A.1) we have

$$2\,\langle w_\alpha^2\, \psi, \mathcal{A}_\varepsilon\, \psi\rangle = 2\|w_\alpha\, \psi\|_{H^{\sigma/2}}^2 - \varepsilon^4 \sum_{z, y} q_\varepsilon(z, y)\, (w_\alpha(z) - w_\alpha(y))^2\, \psi(z)\, \psi(y).$$

Therefore, (A.2) together with $2\,|a\,b|\leqslant a^2+b^2$ give

$$2\,\langle w_\alpha^2\,\psi,\mathcal{A}_\varepsilon\,\psi\rangle\geqslant 2\,\|w_\alpha\,\psi\|_{H^{\sigma/2}}^2-C\,L^{-2}\,\|w_\alpha\,\psi\|_2^2.$$

Taking $L$ sufficiently large proves (A.3).

For the second inequality, we note that

$$\varepsilon^2\sum_z w_\alpha(z)^2\leqslant\varepsilon^2\sum_{z\in\Lambda_\varepsilon}(1+|z|/L)^{-4}\leqslant(1+\varepsilon/L)^4\int_{\mathbb{R}^2}(1+|x|/L)^{-4}\,\mathrm{d}x\leqslant 2\,L^2.$$

Consequently, Cauchy–Schwarz gives

$$\|w_\alpha\,\psi\|_2^4\leqslant\left(\varepsilon^2\sum_z w_\alpha(z)^2\right)\|w_\alpha^{1/2}\,\psi\|_4^4\leqslant 2\,L^2\,\|w_\alpha^{1/2}\,\psi\|_4^4,$$

which proves (A.4). □

**Lemma A.3.** *For every fixed $a,b\geqslant 0$ there is $C_{a,b}<\infty$, independent of $T$, the cutoffs, and the block indices $\alpha,\beta$, such that*

$$\boldsymbol{E}_\beta^\Pi[e^{a\mathcal{T}(\pi)+b\mathcal{N}(\pi)}\,e^{-d_\ell(\Pi_T,\alpha)/L}]\leqslant C_{a,b}\,e^{-d_\ell(\beta,\alpha)/L}\,e^{[b+\kappa(b)]T}. \tag{A.5}$$

*Here $\pi(r)=\Pi_{1+T-r}$ and $\kappa$ is the exponential jump cost* (2.13).

**Proof.** Set $c_*=\sum_{k\geqslant 1}e^{-k/2}<\infty$. For almost every $r\in[0,T]$, let $k(r)$ be fixed by $1+T-r\in I_{k(r)}^T$. In the following, sums over $s<r\leqslant T$ run over the times when the walk $\Pi$ jumps. We have

$$\mathcal{N}(\pi)\leqslant T+1+N_T^\Pi,\qquad\mathcal{T}(\pi)\leqslant c_*+\sum_{0<r\leqslant T}e^{-k(r)/2},$$

where $N_T^\Pi$ is the number of jumps of $\Pi$ during $[0,T]$. Moreover, the triangle inequality and $\Pi_0=\beta$ give

$$d_\ell(\beta,\alpha)\leqslant d_\ell(\Pi_T,\alpha)+\sum_{0<r\leqslant T}d_\ell(\Pi_{r-},\Pi_r).$$

Combining the preceding estimates yields

$$e^{a\mathcal{T}(\pi)+b\mathcal{N}(\pi)}\,e^{-d_\ell(\Pi_T,\alpha)/L}\leqslant e^{ac_*+b(T+1)-d_\ell(\beta,\alpha)/L}\exp\left(\sum_{0<r\leqslant T}F_r\,(\Pi_{r-},\Pi_r)\right), \tag{A.6}$$

where

$$F_r\,(\bar\alpha,\bar\beta):=\frac{d_\ell(\bar\alpha,\bar\beta)}{L}+b+a\,e^{-k(r)/2}.$$

For $0\leqslant s\leqslant T$, set

$$h(s,\bar\alpha):=\boldsymbol{E}_{s,\bar\alpha}^\Pi\exp\left(\sum_{s<r\leqslant T}F_r\,(\Pi_{r-},\Pi_r)\right),$$

where $\boldsymbol{E}_{s,\bar\alpha}^\Pi$ denotes the expectation for the walk started at $\bar\alpha$ at time $s$. By conditioning on the first jump, $h$ satisfies, for almost every $s\in[0,T]$,

$$-\partial_s\,h(s,\bar\alpha)=\sum_{\bar\beta\neq\bar\alpha}j_{\bar\alpha\bar\beta}\,(e^{F_s(\bar\alpha,\bar\beta)}\,h(s,\bar\beta)-h(s,\bar\alpha)),\qquad h(T,\bar\alpha)=1.$$

Set

$$m(r)=\sup_{\bar\alpha}\sum_{\bar\beta\neq\bar\alpha}j_{\bar\alpha\bar\beta}\,(e^{F_r(\bar\alpha,\bar\beta)}-1).$$

The spatially constant function $v(s)=\exp\big(\int_s^T m(r)\,\mathrm{d}r\big)$ is a supersolution of the above backward equation, with the same terminal condition. Since the off-diagonal coefficients are nonnegative, comparison gives $h(0,\beta)\leqslant v(0)$, proving

$$\boldsymbol{E}_\beta^\Pi\exp\left(\sum_{0<r\leqslant T}F_r\,(\Pi_{r-},\Pi_r)\right)\leqslant\exp\left(\int_0^T\sup_{\bar\alpha}\sum_{\bar\beta\neq\bar\alpha}j_{\bar\alpha\bar\beta}\,(e^{F_r(\bar\alpha,\bar\beta)}-1)\,\mathrm{d}r\right). \tag{A.7}$$

By (2.13), (2.6), and $e^{a\rho}-1\leqslant a\,e^a\,\rho$ for $0\leqslant\rho\leqslant 1$,

$$\sum_{\bar{\beta}\neq\bar{\alpha}} j_{\bar{\alpha}\bar{\beta}}\,(e^{F_r(\bar{\alpha},\bar{\beta})}-1)=\sum_{\bar{\beta}\neq\bar{\alpha}} j_{\bar{\alpha}\bar{\beta}}\,(e^{d_\ell(\bar{\alpha},\bar{\beta})/L+b}-1)+e^b\,(e^{a e^{-k(r)/2}}-1)\sum_{\bar{\beta}\neq\bar{\alpha}} j_{\bar{\alpha}\bar{\beta}}\,e^{d_\ell(\bar{\alpha},\bar{\beta})/L}$$
$$\leqslant\kappa(b)+\tfrac{1}{2}a\,e^{a+b}\,e^{-k(r)/2}.$$

We have

$$\int_0^T e^{-k(r)/2}\,\mathrm{d}r\leqslant\sum_{k\geqslant 1}|\{r\colon k(r)=k\}|\,e^{-k/2}\leqslant\sum_{k\geqslant 1}e^{-k/2}\leqslant c_*.$$

Consequently, by (A.7) we arrive at

$$\boldsymbol{E}^{\Pi}_{\beta}\exp\left(\sum_{0<r\leqslant T}F_r\,(\Pi_{r-},\Pi_r)\right)\leqslant\exp\!\big(\kappa(b)\,T+\tfrac{1}{2}a\,e^{a+b}\,c_*\big).$$

Combining this with (A.6), we obtain (A.5) with $C_{a,b}=\exp\,(b+a\,c_*+\tfrac{1}{2}a\,e^{a+b}\,c_*)$. □

# Appendix B Weighted lattice Besov estimates

In this appendix we collect deterministic estimates in weighted Besov spaces on the lattice $\Lambda_\varepsilon$ used to control products involving the singular stochastic terms and the nonlinear remainder.

**Definition B.1.** *Choose a smooth dyadic partition of unity $(\varphi_j)_{j\geqslant -1}$ on $\mathbb{R}^2$ as in [MP19, Section 2.4]. Given the lattice spacing $\varepsilon>0$, let $J_\varepsilon$ be the smallest index such that the support of $\varphi_{J_\varepsilon}$ intersects the boundary of the Brillouin zone and set*

$$\varphi_j^\varepsilon=\varphi_j,\qquad j<J_\varepsilon,\qquad\qquad \varphi_{J_\varepsilon}^\varepsilon=1-\sum_{j<J_\varepsilon}\varphi_j.$$

*Let $K_j^\varepsilon$ be the inverse Fourier transform of $\varphi_j^\varepsilon$ on the infinite lattice. For periodic $f\colon\Lambda_\varepsilon\to\mathbb{R}$, define*

$$\Delta_j^\varepsilon f=K_j^\varepsilon*f,\qquad -1\leqslant j\leqslant J_\varepsilon.$$

*For $\rho\in\mathbb{R}$, $m\geqslant 0$, $q\geqslant 1$, and periodic $F,f\colon\Lambda_\varepsilon\to\mathbb{R}$, define*

$$\|F\|_{\mathcal{C}^\rho(\omega_\alpha^m)}=\sup_{-1\leqslant j\leqslant J_\varepsilon}2^{j\rho}\,\|\omega_\alpha^m\,\Delta_j^\varepsilon F\|_\infty \tag{B.1}$$

*and*

$$\|f\|_{B^\rho_{1,1}(\omega_\alpha^m)}=\sum_{j=-1}^{J_\varepsilon}2^{j\rho}\,\|\omega_\alpha^m\Delta_j^\varepsilon f\|_1,\qquad \|f\|_{B^\rho_{q,\infty}(\omega_\alpha^m)}=\sup_{-1\leqslant j\leqslant J_\varepsilon}2^{j\rho}\,\|\omega_\alpha^m\,\Delta_j^\varepsilon f\|_q \tag{B.2}$$

*We omit the weight in the notation if it is identically equal to 1.*

Note that by Minkowski's inequality and $\omega(z)^m/\omega\,(z-x)^m\leqslant C\,(1+|x|/L)^{m\theta}$, we have

$$\|\omega^m\,(K*f)\|_q\leqslant C\left(\varepsilon^2\sum_{z\in\Lambda_\varepsilon}|K(z)|\,(1+|z|/L)^{m\theta}\right)\|\omega^m\,f\|_q. \tag{B.3}$$

The Fourier symbol of $-\Delta_\varepsilon$ satisfies

$$\frac{4}{\pi^2}|\xi|^2\leqslant\frac{4}{\varepsilon^2}\sum_{i=1}^2\sin^2\!\Big(\frac{\varepsilon\,\xi_i}{2}\Big)\leqslant|\xi|^2$$

on the Brillouin zone. Therefore, for $\rho\in\mathbb{R}$,

$$\|f\|_{H^\rho}^2\asymp\sum_{j=-1}^{J_\varepsilon}2^{2j\rho}\,\|\Delta_j^\varepsilon f\|_2^2. \tag{B.4}$$

**Lemma B.2.** *For $0<\rho<1$, $0\leqslant m\leqslant 3$, and periodic $F, f: \Lambda_\varepsilon \to \mathbb{R}$,*

$$|\langle F, f\rangle| \leqslant C\|F\|_{\mathcal{C}^{-\rho}(\omega_\alpha^m)} \|\omega_\alpha^{-m} f\|_{B_{1,1}^\rho}. \tag{B.5}$$

**Proof.** The result follows by Besov duality from [GH21, Lemmas A.1–A.2]. Note that the periodic weights $\omega_\alpha^m$ may be replaced with smooth periodic weights, without changing the norms. □

**Lemma B.3.** *Let $a, b\in\mathbb{R}$ with $a+b>0$, and $m, n\geqslant 0$ with $m+n\leqslant 3$. For periodic $f, g: \Lambda_\varepsilon\to\mathbb{R}$,*

$$\|f g\|_{\mathcal{C}^{\min\{a,b\}}(\omega_\alpha^{m+n})} \leqslant C\|f\|_{\mathcal{C}^a(\omega_\alpha^m)}\|g\|_{\mathcal{C}^b(\omega_\alpha^n)}. \tag{B.6}$$

**Proof.** Apply the Bony decomposition of the product and use [MP19, Lemma 4.2]. □

**Lemma B.4.** *For every $\alpha\in\mathscr{B}_\ell$, every $\varphi\in\mathcal{H}_\alpha^{\mathrm{loc}}$, and periodic $F, G, \psi: \Lambda_\varepsilon\to\mathbb{R}$,*

$$|\langle F, \varphi^2\rangle| \leqslant C\|F\|_{\mathcal{C}^{-2\eta}(\omega_\alpha^2)}\|\varphi\|_2\|\varphi\|_{H^{2\eta+2\delta}}, \tag{B.7}$$

$$|\langle G\psi, \varphi^2\rangle| \leqslant C\|G\|_{\mathcal{C}^{-\eta}(\omega_\alpha)} \|w_\alpha\psi\|_{H^{\sigma/2}}\|\varphi\|_{H^{\eta+\delta/2}}^2, \tag{B.8}$$

*where $\eta$ is defined by (3.2).*

**Proof.** Write $w=w_\alpha$ and $\omega=\omega_\alpha$. Since $\operatorname{supp}\varphi\subset B_\alpha$, the weights and their inverses are bounded and Lipschitz on a fixed neighbourhood of $\operatorname{supp}\varphi$. Hence, the Sobolev multiplier estimate yields

$$\|\omega^{-1}\varphi\|_{H^b} + \|(w\omega)^{-1/2}\varphi\|_{H^b} \leqslant C\|\varphi\|_{H^b}, \qquad 0\leqslant b\leqslant \sigma/2. \tag{B.9}$$

By (B.5) and [GH21, Lemma A.7],

$$|\langle F, \varphi^2\rangle| \leqslant C\|F\|_{\mathcal{C}^{-2\eta}(\omega^2)}\|(\omega^{-1}\varphi)^2\|_{B_{1,1}^{2\eta}} \leqslant C\|F\|_{\mathcal{C}^{-2\eta}(\omega^2)}\|\varphi\|_2\|\varphi\|_{H^{2\eta+2\delta}}.$$

This proves (B.7). For the second estimate, first note that the paraproduct bounds of [MP19, Lemma 4.2], together with the Sobolev embeddings, yield

$$\|f_1 f_2 f_3\|_{B_{1,1}^\eta} \leqslant C\prod_{i=1}^3 \|f_i\|_{H^{s_i}}, \qquad \eta<s_i\leqslant 1, \quad \sum_{i=1}^3 s_i > 1+\eta.$$

Apply this with $f_1 = w\psi$, $f_2=f_3=(w\omega)^{-1/2}\varphi$ and $s_1=\sigma/2$, $s_2=s_3=\eta+\delta/2$. The bound (B.8) follows then by (B.5) and (B.9). □

**Lemma B.5.** *For every $t\in\mathbb{R}$, $\alpha\in\mathscr{B}_\ell$, and periodic $\psi: \Lambda_\varepsilon\to\mathbb{R}$,*

$$\begin{aligned} &6\,|\langle Z(t), w_\alpha^2\psi^3\rangle| + 6\,|\langle V(t), w_\alpha^2\psi^2\rangle| + 2\,|\langle H(t), w_\alpha^2\psi\rangle| \\ &\quad\leqslant \frac{1}{4}\|w_\alpha\psi\|_{H^{\sigma/2}}^2 + \frac{1}{4}\|w_\alpha^{1/2}\psi\|_4^4 + C\,\mathfrak{A}_\alpha(t)^p. \end{aligned} \tag{B.10}$$

*The constant $C$ is independent of $\varepsilon, \ell, \alpha, t, \psi$, and $\lambda\in(0,1]$.*

**Proof.** Write $w=w_\alpha$ and $\omega=\omega_\alpha$. Since $8\theta<2$, we have

$$\varepsilon^2\sum_{z\in\Lambda_{\varepsilon,\ell}} w(z)^2\,\omega(z)^{-8} \leqslant C.$$

Therefore,

$$\|w^{1/4}\omega^{-1}\|_8 + \|w^{1/2}\omega^{-2}\|_4 + \|w\,\omega^{-3}\|_2 \leqslant C.$$

We use the weighted Besov estimates of [GH21, Appendix A], replacing our periodic weight by comparable smooth periodic weights. Observe also that by (3.2) and the choice of $\delta$,

$$\eta+2\delta<\frac{\sigma}{4}, \qquad 2\eta+2\delta<\frac{\sigma}{2}.$$

By [GH21, Lemmas A.1, A.3 and A.4], we obtain

$$\|w\,\omega^{-1}\psi\|_{H^{\eta+2\delta}} \leqslant C\,\|w^{1/4}\omega^{-1}\|_8\,\|\psi\|_{B_{8/3,\infty}^{\sigma/4}(w^{3/4})} \leqslant C\,\|w\,\psi\|_{H^{\sigma/2}}^{1/2}\,\|w^{1/2}\psi\|_4^{1/2}.$$

Hence, by (B.5) and [GH21, Lemma A.7], we have

$$\begin{aligned}|\langle Z(t), w^2\,\psi^3\rangle| &\leqslant C\,\|Z(t)\|_{\mathcal{C}^{-\eta}(\omega)}\|\,w^{1/2}\,\psi\|_4^2\,\|w\,\omega^{-1}\,\psi\|_{H^{\eta+2\delta}}\\ &\leqslant C\,\|Z(t)\|_{\mathcal{C}^{-\eta}(\omega)}\|\,w\,\psi\|_{H^{\sigma/2}}^{1/2}\|\,w^{1/2}\,\psi\|_4^{5/2}.\end{aligned}$$

Similarly, the quadratic estimate of [GH21, Lemma A.7] gives

$$\begin{aligned}|\langle V(t), w^2\,\psi^2\rangle| &\leqslant C\,\|V(t)\|_{\mathcal{C}^{-2\eta}(\omega^2)}\|\,w\,\omega^{-2}\,\psi\|_2\,\|w\,\psi\|_{H^{2\eta+2\delta}}\\ &\leqslant C\,\|V(t)\|_{\mathcal{C}^{-2\eta}(\omega^2)}\|\,w\,\psi\|_{H^{\sigma/2}}\|\,w^{1/2}\,\psi\|_4.\end{aligned}$$

For the linear term, weighted Besov embedding and Hölder's inequality give

$$\begin{aligned}|\langle H(t), w^2\,\psi\rangle| &\leqslant C\,\|H(t)\|_{\mathcal{C}^{-2\eta}(\omega^3)}\|\,w\,\omega^{-3}\|_2\,\|w\,\psi\|_{H^{2\eta+\delta}}\\ &\leqslant C\,\|H(t)\|_{\mathcal{C}^{-2\eta}(\omega^3)}\|\,w\,\psi\|_{H^{\sigma/2}}.\end{aligned}$$

Combining these estimates and applying Young's inequality, we get

$$\begin{aligned}&6\,|\langle Z(t), w^2\,\psi^3\rangle\,|+6|\langle V(t), w^2\,\psi^2\rangle\,|+2|\langle H(t), w^2\,\psi\rangle|\\ &\quad\leqslant \frac{1}{4}\|w\,\psi\|_{H^{\sigma/2}}^2+\frac{1}{4}\,\|w^{1/2}\,\psi\|_4^4+C\,(\|Z(t)\|_{\mathcal{C}^{-\eta}(\omega)}^8+\|V(t)\|_{\mathcal{C}^{-2\eta}(\omega^2)}^4+\|H(t)\|_{\mathcal{C}^{-2\eta}(\omega^3)}^2).\end{aligned}$$

The result follows now from (3.3), since $p\geqslant 8$. □

# Appendix C Extension to the $O(N)$-invariant model

We briefly discuss the extension to the $\Phi^4$ model with fields valued in $\mathbb{R}^N$, for fixed $N\in\mathbb{N}$. The interaction potential is

$$\frac{\lambda}{4}:|\phi|^4:=\frac{\lambda}{4}\,(|\phi|^4-2\,(N+2)\,C_{\varepsilon,\ell}\,|\phi|^2+N\,(N+2)\,C_{\varepsilon,\ell}^2),\qquad |\phi|^2:=|\phi^1|^2+\cdots+|\phi^N|^2.$$

The corresponding stochastic quantisation equation is

$$\mathrm{d}\Phi_t=-[\mathcal{A}_\varepsilon\,\Phi_t+\lambda\,(|\Phi_t|^2-(N+2)\,C_{\varepsilon,\ell})\,\Phi_t]\,\mathrm{d}t+\sqrt{2}\,\mathrm{d}W_t,\qquad \Phi=(\Phi^1,\ldots,\Phi^N),$$

where the components of $W=(W^1,\ldots,W^N)$ are independent.

Let us explain how to adapt the construction of the stochastic objects from Section 3. We define $X$ to be the stationary vector-valued Gaussian convolution and set

$$(\partial_t+\mathcal{A}_\varepsilon)\,Y=\,:|X|^2\,X:\,=(|X|^2-(N+2)\,C_{\varepsilon,\ell})\,X,\qquad Z=X-\lambda\,Y.$$

The scalar Wick square $X^{:2:}$ is replaced by the matrix $B_X=X\otimes X-C_{\varepsilon,\ell}\,I_N$. Summing up, the required enhancement consists of

$$X^k,\qquad (B_X)^{ij},\qquad Y^k,\qquad (B_X)^{ij}\,Y^k,\quad 1\leqslant i,j,k\leqslant N,$$

with the same weighted Besov regularities as $X$, $X^{:2:}$, $Y$, and $X^{:2:}Y$ in the scalar case.

To see why the deterministic estimates have the same structure, put $\Phi^u=Z+\Psi^u$ and $B_Z=Z\otimes Z-C_{\varepsilon,\ell}\,I_N$, where $I_N$ is the identity matrix. Expansion of the original cubic nonlinearity produces the coercive term $\lambda|\Psi^u|^2\,\Psi^u$, quadratic in $\Psi^u$ term involving $Z$, a linear term involving $B_Z$, and a forcing coinciding with a linear combination of terms of the types $B_X Y$, $XY^2$, and $Y^3$. In order to obtain a desired energy bound we can apply componentwise the product estimates in Lemma B.5.

The main point in the local Rayleigh estimate is that the Hessian of the interaction is matrix-valued but has the decomposition

$$\begin{aligned}&(|\Phi^u|^2-(N+2)\,C_{\varepsilon,\ell})\,I_N+2\,\Phi^u\otimes\Phi^u\\ &\quad=|\Psi^u|^2\,I_N+2\,\Psi^u\otimes\Psi^u+(\mathrm{tr}B_Z)\,I_N+2\,B_Z+2\,(Z\cdot\Psi^u)\,I_N+2\,(Z\otimes\Psi^u+\Psi^u\otimes Z),\end{aligned}$$

where the first two terms have quadratic form $|\Psi^u|^2|v|^2 + 2\,(\Psi^u\cdot v)^2 \geqslant 0$. The remaining terms are controlled componentwise by the quadratic and trilinear bounds of Lemma B.4. This gives the analogue of Proposition 4.3, with the vector-valued local Rayleigh quotient defined using the symmetric quadratic form above.

Finally, the localisation and block-transfer argument acts on the scalar energies $e_\alpha = \|\chi_\alpha v\|_2^2$ and requires no positivity of a matrix-valued propagator. The block walk, Boué–Dupuis argument (with vector-valued controls), and covariance representation do not require any modifications. Thus, the proofs of Theorems 1.1 and 1.2 require only very minor modifications.

## Appendix D Synchronisation by noise

Let us explain how to conclude almost sure exponential synchronisation of solutions of the infinite volume continuum equation driven by the same noise. This requires a minor extension of Theorem 1.2. Recall that $\mathfrak{S}_1 = \{\Phi_1^\phi : \phi \in L^2\}$. We claim that the estimate (1.7) is also valid when the supremum over $\psi \in \mathfrak{S}_1$ is replaced with the supremum over its closed convex hull $\bar{\mathrm{co}}(\mathfrak{S}_1)$.

Let us first show that (4.4) holds uniformly over $\psi \in \bar{\mathrm{co}}(\mathfrak{S}_1)$. By convexity of the squared $L^2$ norm,

$$\sup_{\psi \in \bar{\mathrm{co}}(\mathfrak{S}_1^u)} \|w_\alpha(\psi - Z_1)\|_2^2 \leqslant \sup_{\phi \in L^2} \|w_\alpha(\Phi_1^{\phi,u} - Z_1)\|_2^2,$$

where $\mathfrak{S}_1^u := \{\Phi_1^{\phi,u} : \phi \in L^2\}$ and $Z$ is the process introduced in Section 3. The RHS is bounded by the coming-down estimate (4.7). Consequently, applying the energy inequality (4.2) to solutions restarted at time one with the initial condition $\psi \in \bar{\mathrm{co}}(\mathfrak{S}_1^u)$ yields (4.4). Thus, the proof of Lemma 5.1 applies uniformly over $\bar{\mathrm{co}}(\mathfrak{S}_1^u)$. Under a Cameron–Martin shift, $\mathfrak{S}_1$ is transformed into $\mathfrak{S}_1^u$, so Proposition 5.2 and the argument of Proposition 2.2 prove (1.7) with the supremum over $\psi \in \bar{\mathrm{co}}(\mathfrak{S}_1)$. This proves the convex hull extension of Theorem 1.2.

To establish synchronisation, we study the random diameter of the solution set at time $t$,

$$\mathfrak{D}_\alpha(t) := \sup_{\phi_1, \phi_2} \|\chi_\alpha\,(\Phi_t^{\phi_1} - \Phi_t^{\phi_2})\|_{L^2},$$

where the supremum runs over all initial conditions. For arbitrary initial conditions $\phi_1, \phi_2$, write $\psi_i = \Phi_1^{\phi_i}$. Differentiability of the solution map gives

$$\Phi_{1+T}^{\phi_1} - \Phi_{1+T}^{\phi_2} = \int_0^1 U_{1,1+T}^{(1-a)\,\psi_2 + a\,\psi_1}(\psi_1 - \psi_2)\,\mathrm{d}a.$$

Since $(1-a)\,\psi_2 + a\,\psi_1 \in \bar{\mathrm{co}}(\mathfrak{S}_1)$, we obtain

$$\mathfrak{D}_\beta(1+T) \leqslant \sum_{\alpha \in \mathscr{B}_\ell} \mathcal{K}_{\beta\alpha}(T)\,\mathfrak{D}_\alpha(1), \qquad \mathcal{K}_{\beta\alpha}(T) := \sup_{\psi \in \bar{\mathrm{co}}(\mathfrak{S}_1)} \|\chi_\beta U_{1,1+T}^\psi \chi_\alpha\|_{L^2 \to L^2}.$$

By (4.7), (2.1), and Proposition 3.3,

$$\sup_{\varepsilon,\ell,\alpha} \mathbb{E}\,\mathfrak{D}_\alpha(1)^2 \leqslant C_\lambda.$$

Consequently, by the convex-hull extension of (1.7), discussed above, and the Cauchy–Schwarz inequality, we obtain

$$\mathbb{E}\,\mathfrak{D}_\beta(1+T) \leqslant C_\lambda\,e^{-T/8} \sum_{\alpha \in \mathscr{B}_\ell} e^{-d_\ell(\beta,\alpha)/(2L)} \leqslant C_\lambda\,e^{-T/8}.$$

This implies almost sure exponential contraction of the diameter with a deterministic decay rate independent of the cutoffs. Unlike [GT20; DHYZ25], this argument does not use order preservation and is therefore applicable to the vector-valued model.